\documentclass[12pt]{amsart}
\usepackage{amssymb, mathrsfs, color}

\numberwithin{equation}{section}

\newtheorem{thm}{Theorem}[section]
\newtheorem{cor}[thm]{Corollary}
\newtheorem{lem}[thm]{Lemma}
\newtheorem{prop}[thm]{Proposition}
\theoremstyle{remark}
\newtheorem{ex}[thm]{Example}

\newtheorem{ques}[thm]{Question}
\newtheorem{rk}[thm]{Remark}

\newcommand{\HH}{\mathcal{H}}

\newcommand{\cH}{{\mathcal{H}}}
\newcommand{\C}{\mathbb{C}}

\newcommand{\D}{\mathbb{D}}
\newcommand{\T}{\mathbb{T}}
\newcommand{\intt}{\int_{\mathbb{T}}}

\newcommand{\ov}{\overline}
\newcommand{\onb}{\textsc{onb}}
\newcommand{\fr}{\frac}
\newcommand{\lam}{\lambda}
\newcommand{\va}{\varphi}
\newcommand{\al}{\alpha}

\newcommand{\ran}{\textup{ran}}
\newcommand{\Cla}{C_{\lam,u}}
\newcommand{\Clae}{C_{\lam,u,E}}
\newcommand{\uet}{\textsc{uet}}

\begin{document}

\title[Complex symmetry of Toeplitz operators]{Complex symmetry of Toeplitz operators}

\author[Y. Chen]{Yong Chen}
\address{Department of Mathematics, Hangzhou Normal University, Hangzhou, 311121, P.R.China}
\email{ychen@hznu.edu.cn}

\author[Y. J. Lee]{Young Joo Lee}
\address{Department of Mathematics, Chonnam National University, Gwangju, 61186, Korea}
\email {leeyj@chonnam.ac.kr}

\author[Y. Zhao]{Yile Zhao}
\address{Department of Mathematics, Hangzhou Normal University, Hangzhou, 311121, P.R.China}
\email{yilezhao@hznu.edu.cn}

\thanks{}
\date{\today}

\begin{abstract} 
We introduce a new class of conjugations and characterize
complex symmetric Toeplitz operators  on the Hardy space with respect  to those conjugations. Also, we prove  that complex symmetricity and \uet~ property are the same for a certain class of Toeplitz operators.
We also discuss the
analytic symmetricity for Toeplitz operators.
Our results extend several known results by providing unified ways of treating them.
\end{abstract}

\maketitle



\section{Introduction}

Given  a complex separable Hilbert space $\cH$ endowed with the inner product $\langle\ ,\ \rangle_\HH$,
a conjugate-linear map $C: \cH\rightarrow\cH$ is called a conjugation if $C^2=I$ and $\langle Cx,Cy\rangle_\HH=\langle y, x\rangle_\HH$ for all $x, y\in \HH$.
Given a bounded linear operator $T$ and conjugation $C$ on $\HH$, $T$ is called  $C$-symmetric if $CTC=T^*$. Also, we say that $T$ is {\it complex symmetric} if $T$ is  $C$-symmetric for some conjugation $C$ on $\HH$. The term ``complex symmetric" is due to the
fact that $T$ is complex symmetric if and only if $T$ is unitarily equivalent to
a symmetric matrix with complex entries regarded as an operator acting on an $\ell^2$-space
of the appropriate dimension.

The class of complex symmetric operators is very large and includes not only normal operators but also lots of nonnormal operators containing  Hankel operators, truncated Toeplitz operators and Volterra integration operators.
In order to obtain a better understanding for complex symmetric operators, a lot of effort are devoted to the study of complex symmetricity for several kinds of operators containing  partial isometries \cite{SG6}, weighted shifts \cite{ZhuLi},  triangular operators \cite{ZhuOAM} and composition operators \cite{GaoZhou,JKKL,LimK,Thompson,WangHan}.

Recently,  it has been studied the problem of when a Toeplitz operator on the Hardy space is complex symmetric. We first recall the Hardy space.
Let $\T$ be the boundary of the open unit disk $\D$ in the complex plane $\C$. We let $L^2=L^2(\T,\sigma)$ be the usual Lebesgue space on $\T$ where $\sigma$ is the normalized Haar measure on $\T$.
The classical Hardy space $H^2$ is consisting of all analytic functions $f$ on $\D$ satisfying
$$
\sup_{0\le r<1}\int_\T |f(r\xi)|^2d\sigma(\xi)<\infty.
$$ 
As is well known, the Hardy space $H^2$ is isometrically identified with a closed subspace of $L^2$ via the boundary functions. Thus, we will use the same letter for a function in $H^2$ and its boundary function in $L^2$.
In particular,  $H^2$ is a Hilbert space with the inner product
$$
\langle f, g\rangle
=\int_\T f\bar g\, d\sigma
$$
for functions $f,g \in H^2$.
Also, it is well known that the Hilbert space orthogonal projection $P$ from $L^2$ onto $H^2$ is realized by the integral operator
$$
P\psi(z)=\langle \psi,K_z\rangle=\intt \fr{\psi(\xi)}{1-z\bar\xi}\,d\sigma(\xi),\qquad z\in \D
$$
for $\psi\in L^2$ where $K_z$ is the  reproducing kernel
on $H^2$ given by
$$
K_z(\xi)=\frac1{1-\bar z\xi},\qquad \xi\in \T.
$$
Given $\phi\in L^\infty(\T)$, the {\it Toeplitz operator} $T_\phi$ with symbol $\phi$ is defined  by
$$
T_\phi f =P(\phi f)  
$$
for $f\in H^2.$ Clearly the Toeplitz operator $T_\phi$ is a bounded linear operator on $H^2$.

As is well known, for $\phi\in L^\infty(\T)$, 
$T_\phi$ is normal on $H^2$ if and only if $\phi=f+\ov{\gamma f}$ for some $f\in H^\infty$ and $\gamma\in\T$; see \cite{BH} for details. So, such $\phi$ induces a complex symmetric Toeplitz operator.
Guo and Zhu \cite{GZ} characterized a certain trigonometric symbol for which the corresponding Toeplitz operator is complex symmetric.
They also raised an interesting question at the same paper: {\it Characterize complex symmetric Toeplitz operators  on $H^2$.}
Generally it seems to be difficult to describe when a Toeplitz operator is complex symmetric.
Recently several kinds of conjugations  have been considered and then
 the characterizing problem of
complex symmetric Toeplitz operators with respect to those conjugations has been studied.

Recently, Ko and Lee \cite{KL} considered a class of conjugations and studied complex symmetric Toeplitz operators.
More explicitly, given $\lam\in\T$,  define
$$
C_\lam h(z)=\ov{h(\lam\bar z)},  \qquad z\in\T
$$ 
for functions $h\in L^2$. 
In other words, if we write $h=\sum_{j=-\infty}^\infty c_n z^n$ for the Fourier series expansion  of $h$, $C_\lambda h$ can written as
$$
C_\lambda h=\sum_{n=-\infty}^\infty \overline{c_n} \overline{\lambda^n}z^n.
$$
Then, one can easily see that $C_\lam$ is a conjugation on both $L^2$ and $H^2$. Given $\phi\in L^\infty(\T)$ and $\lambda\in \T$,  Ko and Lee proved that the Toeplitz operator $T_\phi$ is $C_\lam$-symmetric on $H^2$ if and only if $\phi=f+\ov{C_\lam f}$ for some $f\in H^2$, which is also equivalent to $\phi(\lambda z)=\phi(\bar z)$ for all $z\in \D$; see \cite{KL} for details.

Recently, a more general class of conjugations  induced by sequences has been introduced in \cite{LiYangLu}. Given a sequence $\alpha=(\alpha_0, \alpha_1,\alpha_2, \cdots)$ where $\alpha_n\in\T$, define a conjugation $C_\alpha$ by 
$$
C_\alpha h=\sum_{n=0}^\infty \overline{c_n}\alpha_n  z^n
$$
for $h=\sum_{n=0}^\infty c_n z^n\in H^2.$
Specially, $C_\alpha=C_\lambda$ when $\alpha_n=\overline{\lambda^n}$ for each $n$.
Such a conjugation is known to be a canonical one. Recall that
a conjugation $C$ is called a  canonical conjugation if  $Cz^n\in \vee\{z^n\}$ for all $n=0,1,2,\cdots$ where $\vee$ denotes the closed linear span. 
Given a conjugation $C$ on $H^2$, we note that $C$ is a canonical conjugation if and only if $C=C_\alpha$ for some $\alpha=\{\alpha_n\}$ where $\alpha_n\in\T$.
An operator  is called  {\it canonically symmetric} if it is complex symmetric with respect to a canonical conjugation.

In a recent paper \cite{LiYangLu}, the authors characterized canonically symmetric Toeplitz operators in terms of the notion of  geometric pair sequences; see Theorem 2.11 of \cite{LiYangLu} for details.
Very recently, it turned out that
canonically symmetricity and 
$C_\lambda$-symmetricity are the same for Toeplitz operators.
More explicitly, given $\phi\in L^\infty(\T)$,
  it is proved that  $T_\phi$ is canonically symmetric  if and only if  $T_\phi$ is $C_\lam$-symmetric for some $\lambda\in \T$, which provides a finer version of the result in \cite{LiYangLu}; see \cite{BCZ} for details and related results.

Motivated by these results,  we introduce a more general class of conjugations induced by a sequence and an inner function. 
To be more precise, we first need some notations.
Let $\mathscr M$ be the set of all sequences $\alpha=(\alpha_0, \alpha_1,\alpha_2, \cdots)$ where each $\alpha_n\in\T$. Also, the notation $H^\infty$ denotes the space of all bounded analytic functions on $\D$ and a function $u\in H^\infty$ is said to be inner if  $|u|=1$ a.e. on $\T$.

Let $\alpha=\{\alpha_m\}\in\mathscr M$ and $u$ be a nonconstant inner function. Choose an orthonormal basis (briefly, \onb) $E=\{e_j: j\in\Lambda \}$ of the model space $H^2\ominus u H^2$. By the Wold decomposition theorem,
we have
$$
H^2=\vee\{ u^m e_j:  j\in \Lambda,\, m=0,1,2,\cdots\}.
$$
We define $C_{\alpha,u,E}$  on $H^2$ by 
\begin{equation*}
C_{\alpha,u,E}\Bigl( \sum_{m, j} \mu_{m, j} u^m e_j\Bigr)=\sum_{m, j} \overline{\mu_{m, j}}
\alpha_m u^m e_j.
\end{equation*}
Then, one can check that $C_{\alpha,u,E}$ is a conjugation on $H^2$.
Note that the $C_{\alpha,u,E}$-symmetricity of an operator depends on the choice of an \onb~ $E$; see Example \ref{p-7}.
We say that an operator $T$ on $H^2$ is $C_{\alpha,u}$-{\it symmetric} if $T$ is $C_{\alpha,u,E}$-symmetric for every \onb~ $E$ of $H^2\ominus uH^2$.

In the special case when $u=z$ and $E=\{1\}$, the corresponding conjugation 
$C_{\alpha,u,E}$ is a canonical conjugation $C_\alpha$.
More specially, if $\alpha$ is given by a geometric sequence such as  $\alpha_n=\bar\lam^n$ for  $\lam\in\T$, the corresponding conjugation $C_{\alpha,u,E}$ is just conjugation $C_\lam$.
The $C_{\alpha,u}$  will be denoted by $C_{\lam,u}$  if $\alpha$ is given by a geometric sequence induced by  $\lam\in\T$ above.

For conjugations $C_{\alpha,u,E}$, we first obtain the following characterization of when a Toeplitz operator is $C_{\alpha,u}$-symmetric on $H^2$.

\begin{thm}
\label{C1-1}
Let $\phi\in L^\infty(\T)$ and $u$ be inner with a zero in $\D$. Then the following statements are equivalent:
\begin{enumerate}
\item[(a)] $T_\phi$ is $C_{\alpha, u}$-symmetric for some $\alpha\in\mathscr M$.
\item[(b)] $T_\phi$ is $\Cla$-symmetric for some $\lam\in\T$.
\item[(c)] $\phi=f\circ u+\ov{C_\lam(f)\circ u}$ for some $\lam\in\T$ and $f\in H^2$.
\end{enumerate}
\end{thm}

We were not able to prove the same in the case when the inner function $u$ has no zero.
The proof of Theorem \ref{C1-1} will be given at Section \ref {S:CLambda}. 
As a special case of  when an inner function is given by a single M\"obius transform, 
 we extend and recover several known results  in \cite{BCZ} and \cite{KL}.

Recall that a bounded linear operator $T$ on a Hilbert space $\HH$ is called \uet\ if $T$ is unitary equivalent to $CT^*C$ for some conjugation $C$ on $\HH$. The notion of \uet\  was inspired by Halmos \cite{Halmos} and has been studied extensively in \cite{GJZ} and \cite{GZ}. In particular, it plays an important role in describing the block structure of complex symmetric operators.

Note that each complex symmetric operator is \uet ~ but the converse is not true generally; see \cite{GT} or \cite{GJZ} for details.
However, it turns out that the converse is also true for a certain class of Toeplitz operators on $H^2$. In \cite{BCZ}, it was recently shown that, for a Toeplitz operator with certain trigonometric symbol,  the complex symmetricity and \uet~ 
 property are the same, which is equivalent to that
it is either normal or $C_\lam$-symmetric for some $\lambda\in \T$. Having this result,  the authors raised the following questions.

\begin{ques}
\label{C:main}
(a) Is a \uet~ Toeplitz operator complex symmetric?

(b) Is a non-normal complex symmetric Toeplitz operator $C_\lam$-symmetric for some $\lambda\in \T$?
\end{ques}

For Question \ref{C:main},  we consider a more general class of Toeplitz operators with 
non-trigonometric  symbols.
For a point $a\in\D$, we let
$$
\va_a(z)=\fr{a-z}{1-\bar az},\qquad z\in\D
$$ 
denote the usual M\"obius transformation of $\D$. Also, a function of the form
$B=\gamma \va_{a_1}\va_{a_2}\cdots\va_{a_N}$  where $\gamma \in\T$ and $a_j\in\D$ is called
a  finite Blaschke product.

Our next result below provides another evidence for a positive answer to (a) of Question \ref{C:main} and a negative answer to (b) of Question \ref{C:main}.

\begin{thm}\label{T:main-2}
Let $\phi=B_1+\ov{\mu B_2\circ B_1}$ where $\mu\in\C$ and $B_1, B_2$ are finite Blaschke products. Then the following assertions  are equivalent:
\begin{enumerate}
\item[(a)] $T_\phi$ is complex symmetric.
\item[(b)] $T_\phi$ is \uet.
\item[(c)] $\mu B_2=q\varphi_b$ for some $q\in\T$ and $b\in \D$ satisfying $qb=\bar b$.
\end{enumerate}
\end{thm}

Recently, the complex symmetricity of Toeplitz operators has been studied in terms of a certain geometric property of the symbol. 
Recall that a closed curve $\gamma$ in $\mathbb C$ is called nowhere winding if its winding number at  $z$ is 0 for all
$z\in\mathbb C\setminus \gamma$; see Section \ref{S:Blaschke} for details.
Noor  \cite{Noor}  recently proved that
the image of $\T$ under a continuous symbol of a complex symmetric Toeplitz operator is nowhere winding. Based on this result, Noor posed a question of whether the converse is also true.

Thus, in view of Question \ref{C:main} together with Noor's question above, we notice that there might be some connections between complex symmetricity,  \uet~ condition and nowhere winding property.
We show that the image
$\phi(\T)$ of a \uet~ Toeplitz operator with continuous symbol $\phi$ is nowhere winding either; see Proposition \ref{p-4}. Also, for symbols of the form $\phi=z+\ov{q\va_b}$ with $|q|\le 1$
and $b\in\D$, we show that $T_\phi$ is \uet\ if and only if $T_\phi$ is complex symmetric, which is also equivalent to
that $\phi(\T)$ is nowhere winding; see Proposition \ref{th-4}. This gives a positive example to Noor's
question mentioned above. 
Those results including Theorem \ref{T:main-2} will be proved 
at Section \ref{S:Blaschke}.

Finally, we have another observation for complex symmetric Toeplitz operators.
For $f,g\in H^\infty$, put $\phi=f+\bar g$. If $T_\phi$ is complex symmetric
for some conjugation $C$ on $H^2$, then
$$
CT_fC+CT_{\bar g}C=T_{\bar f}+T_g.
$$
Since $T_g,CT_fC$ are hyponormal and $CT_{\bar g}C, T_{\bar f}$ are co-hyponormal, we naturally
guess that $CT_fC$ is a Toeplitz operator with  analytic symbol and $CT_fC=T_g$. But it is not the case  generally; see Example \ref{E:counterEg}.
So we need a definition.

Given $\phi\in L^\infty(\T)$,  we say that $T_\phi$ is {\it analytically symmetric} if there exist $f,g\in H^\infty$ and a conjugation $C$ on $H^2$ such that $CT_fC=T_g$ and $\phi=f+\bar g$.
In this case, since $CT_gC=T_f$ and $CT_{\bar{g}}C=T_{\bar{f}}$, we see $$CT_{\phi}C=T_{g+\bar{f}}=T^*_{\phi}.$$
Thus, analytic symmetric Toeplitz operators are complex symmetric. But the converse is not true in general; see  also Example \ref{E:counterEg}.
So, we have  a  question below naturally.

\begin{ques}
\label{ansym}
For which $\phi\in L^\infty(\T)$, is the Toeplitz operator $T_\phi$  analytically symmetric on $H^2$?
\end{ques}

In Section \ref{S:Regular}, we study Question \ref{ansym} and first show that,
for a Toeplitz operator with general
trigonometric symbol,
the analytic symmetric property is equivalent to the  $C_\lambda$-symmetry for some $\lambda\in \T$; see Theorem \ref{p-6}. 
Also, we characterize
an analytic symmetric Toeplitz operator when the analytic part of the symbol is a polynomial of degree 2; see Theorem \ref{p-new-1}.
Finally, in conjunction with the notion of  analytic symmetric Toeplitz operator, we characterize conjugations
$C$ on $H^2$ for which $CT_fC$ is a Toeplitz operator with an analytic symbol for any $f\in H^\infty$; see Theorem \ref{p-8}.


\section{$C_{\alpha,u}$-symmetric Toeplitz operators}
\label{S:CLambda}

In this section we prove Theorem \ref{C1-1} and
obtain immediate  consequences in case when the inner function is given by a single M\"obius transform.

Given $\alpha=\{\alpha_m\}\in\mathscr M$ and an \onb~ $E=\{e_j: j\in\Lambda \}$ of $H^2\ominus u H^2$ corresponding to a nonconstant inner function $u$, we recall
the conjugation $C_{\alpha,u,E}$  on $H^2$ defined by
\begin{equation*}
\label{e-e2}
C_{\alpha,u,E}: \  u^m e_j\mapsto  \alpha_m u^m e_j, \qquad   j\in \Lambda,\, m=0,1,2\cdots.
\end{equation*}
We also recall that $C_{\alpha,u}$ is denoted by $C_{\lam,u}$  if $\alpha=\{{\bar\lam}^m\}$ for some $\lambda\in\T$.

We first observe basic results which will be useful in our characterizations. One of them shows
that $T_\phi$ is $C_\lam$-symmetric if and only if $T_{\phi\circ u}$ is $\Cla$-symmetric,
which  not only provides a much larger class of complex symmetric Toeplitz operators but  also reminds that (b) of Question \ref{C:main} may be true only for Toeplitz operators with trigonometric symbols.

\begin{prop}
\label{C1}
Let $\phi\in L^\infty(\T)$ and $u$ be a non-constant inner function. Then the following statements are equivalent:
\begin{itemize}
\item[(a)] $T_\phi$ is canonically symmetric.
\item[(b)] $T_\phi$ is $C_\lam$-symmetric for some $\lam\in\T$.
\item[(c)] $T_{\phi\circ u}$ is $\Cla$-symmetric for some $\lam\in\T$.
\item[(d)] $T_{\phi\circ u}$ is $\Clae$-symmetric for some $\lam\in\T$ and \onb~ $E$ of $H^2\ominus uH^2$.
\item[(e)] $\phi=f+\ov{C_\lam f}$ for some $\lam\in\T$ and $f\in H^2$.
\item[(f)] $\phi(\lambda z)=\phi(\overline{z})$ for some $\lam\in\T$.
\end{itemize}
\end{prop}

\begin{proof} First, implications 
(a)$\Longleftrightarrow$(b)$\Longleftrightarrow$(e)$\Longleftrightarrow$(f) follow from \cite{BCZ} and \cite{KL} as mentioned before. Also, (c)$\Longrightarrow$ (d) is obvious. Thus, it remains to prove implications (e) $\Longrightarrow$ (c) and (d) $\Longrightarrow$ (f).

First assume (e) and show (c).
Fix an \onb~ $E=\{e_j:j \in\Lambda\}$ of $H^2\ominus uH^2$.
Write $\phi
=\sum_{n=-\infty}^\infty a_n z^n$ for the Fourier series expansion  of $\phi$.
Note $P(u^\ell e_j)=0$ for all integer $\ell<0$ and $j\in\Lambda$. It follows that
\begin{align*}
\Clae T_{\phi\circ u}\Clae (u^me_j)
&=\Clae T_{\phi\circ u}( \bar\lam^m u^me_j)\\
&=\lam^m \Clae 
P\Bigl(\sum_{n=-\infty}^\infty a_n u^{n+m}e_j\Bigr)\\
&=\lam^m \Clae 
\Bigl(\sum_{n=-m}^\infty a_n u^{n+m}e_j\Bigr)\\
&=\lam^m
\sum_{n=-m}^\infty \ov{a_n}\ov{\lambda^{n+m}} u^{n+m}e_j\\
&=
P\Bigl(\sum_{n=-\infty}^\infty \ov{a_n}\ov{\lambda^{n}} u^{n+m}e_j\Bigr)\\
&=
P\bigl[ ((C_\lambda\phi)\circ u) u^me_j\bigr]\\
&=
T_{(C_\lambda\phi)\circ u}
 (u^me_j)
\end{align*}
for all $m, j$ and hence
\begin{align}
\label{e-new-ee2}
\Clae T_{\phi\circ u}\Clae 
=T_{(C_\lambda\phi)\circ u}.
\end{align}
 Now, since $\phi=f+\ov{C_\lam f}$ by the assumption, 
 (\ref {e-new-ee2}) shows that
\begin{align*}
\Clae T_{\phi\circ u}\Clae
&=T_{C_\lam(f+\overline{C_\lambda f})\circ u}\\
&=T_{C_\lam(f+C_\lambda \bar f)\circ u}\\
&=T_{(C_\lam f+\bar f)\circ u}\\
&=T^*_{\phi\circ u}.
\end{align*}
So $T_{\phi\circ u}$ is $\Clae$-symmetric for any \onb~ $E$ and  hence (c) follows.

Now assume (d) and show (e).
Since $T_{\phi\circ u}$ is $C_{\lam,u,E}$-symmetric, it follows from (\ref{e-new-ee2}) that
$
T_{(C_\lam \phi)\circ u}=T_{\ov{\phi\circ u}},
$
which implies $(C_\lam \phi)\circ u=\ov{\phi\circ u}$. Hence $C_\lam \phi=\ov\phi$ and (f) holds. The proof is complete.
\end{proof}


Before we prove Theorem \ref{C1-1}, we need some facts.
Fix a point $a\in \D$. Let $e_a=\sqrt{1-|a|^2}{K_a}$ be the normalized kernel and define an operator $U_a$ on $H^2$ by $U_af=(f\circ\va_a) e_a$ for $f\in H^2$.
Then, using the fact $(e_a\circ\va_a) e_a=1$, one can easily check that $U_a$ is a unitary operator and
\[U_a^*=U_a^{-1}=U_a, \qquad U_aT_\phi U_a=T_{\phi\circ\va_a}.\]
Thus $T_\phi\cong T_{\phi\circ\va_a}$. Here and in what follows, the notation $\cong$ denotes unitary equivalence as usual. Also,
for a conjugation $C$ on $H^2$, noting
 $U_aCU_a$ is another conjugation on $H^2$,
we see that  $T_\phi$ is $C$-symmetric if and only if $T_{\phi\circ\va_a}$ is  $U_aCU_a$-symmetric.

Now we are ready to prove Theorem \ref{C1-1}.

\begin{proof}[Proof of Theorem \ref{C1-1}]
First assume (c) and show (b). 
By an application of (\ref{e-new-ee2}), we get
$
\Clae T_\phi\Clae=T^*_\phi
$
for any \onb~ $E$ of $H^2\ominus uH^2$, so (b) holds.
Also, implication (b) $\Longrightarrow$ (a) follows from taking $\alpha_m=\overline{\lambda^m}$ for each $m$.

We finally assume (a) and show (c). 
Since $u$ has a zero in $\D$, let $u(a)=0$ for some $a\in \D$.
It is easy to check that $E$ is an \onb~ of $H^2\ominus uH^2$
if and only if $U_aE$ is an \onb~ of $H^2\ominus (u\circ\va_a) H^2$. Also, we note that
$U_a C_{\alpha,u,E} U_a=C_{\alpha, u\circ\va_a, U_aE}$. Thus, as mentioned before, $T_\phi$ is $C_{\alpha,u}$-symmetric if and only if $T_{\phi\circ\va_a}$ is $C_{\alpha, u\circ\va_a}$-symmetric. Since $u\circ\va_a(0)=0$, the observation above shows that we may assume $a=0$ without loss of generality. Thus, $u(0)=0$ and so $1\in H^2\ominus uH^2$. Choose an \onb~ $E=\{1, e_1, e_2,\cdots\}$ of $H^2\ominus uH^2$.
Since $H^2=\bigvee_{m\geq 0} u^mE$, we may write
$$
\phi=\sum_{j\geq 0}(a_ju^j+b_j\bar u^j)+\sum_{k\geq 1}\sum_{j\geq 0}(\eta_{k,j}u^je_k+\beta_{k,j}\bar u^j\ov{e_k}).
$$
Let $\alpha=\{\alpha_m\}.$
Since $T_\phi$ is $C_{\alpha,u,E}$-symmetric,  we have for any integers $n, m\ge 0$,
$$
\langle C_{\alpha,u,E} T_\phi C_{\alpha,u,E} u^m,u^ne_\ell\rangle=\langle T_{\bar\phi}u^m,u^n e_\ell\rangle
$$ 
for all integers $\ell\ge 1$, which implies that for $n\geq m$,
$$
\ov{\alpha_m}\alpha_n\ov{\eta_{\ell,n-m}}=\ov{\beta_{\ell,n-m}}, \qquad\ell\geq 1.
$$
It follows that
\begin{equation}
\label{e-21}
\beta_{\ell,j}=\alpha_m\ov{\alpha_{m+j}}\,\eta_{\ell,j}
\end{equation}
for all $m\ge 0$, $\ell\geq 1$ and $j\geq 0$.

Fix $\gamma\in\T$ and set $E_\gamma=\{1, \gamma e_1, \gamma e_2,\cdots\}$. It is clear that $E_\gamma$
is also an \onb~ of $H^2\ominus uH^2$. We can rewrite $\phi$ with respect  to $E_\gamma$ as
$$
\phi=\sum_{j\geq 0}(a_ju^j+b_j\bar u^j)
+\sum_{k\geq 1}\sum_{j\geq 0}\Big(\bar\gamma\eta_{k,j} u^j\gamma e_k+\gamma\beta_{k,j}\bar u^j\ov{\gamma e_k}\Big).
$$
Note $T_\phi$ is also $C_{\alpha,u,E_{\gamma}}$-symmetric. 
Using the similar arguments as we have done above, we obtain that
$
\gamma\beta_{\ell,j}=\alpha_m\ov{\alpha_{m+j}}\bar\gamma\eta_{\ell,j}
$
and hence
$$
\beta_{\ell,j}=\bar\gamma^2\alpha_m\ov{\alpha_{m+j}}\,\eta_{\ell,j}
$$
for all $m\ge 0$, $\ell\geq 1$ and $j\geq 0$.
Noting the above holds for any $\gamma\in\T$, we have
$
\beta_{\ell,j}=\eta_{\ell,j}=0
$
for all $\ell\geq 1$ and $j\geq 0$.
Hence we have
$$
\phi=\sum_{j=0}^\infty(a_ju^j+b_j\bar u^j).
$$
Also, using the equality
$$
\langle C_{\alpha,u,E} T_\phi C_{\alpha,u,E} u^m,u^n\rangle=\langle T_{\bar\phi}u^m,u^n\rangle,
$$
we have
$
\ov{\alpha_m}\alpha_n\ov{a_{n-m}}=\ov{b_{n-m}}
$
for all $n>m\geq 0$, which implies
$$
b_j=\alpha_m\ov{\alpha_{m+j}}\, a_j,\qquad j\geq 1, m\geq 0.
$$
Since $\alpha\in\mathscr M$, we have $|b_j|=|a_j|$ for any $j\ge 1$ and then
\begin{equation}\label{e2}
\fr{\al_{m+j}}{\al_m}=\fr{a_j}{b_j}
\end{equation}
for all $m\geq 0$ and $j\ge 1$ with $a_j\not=0.$
If $a_j=0$ for all $j\ge 1$, then clearly (c) holds. So, suppose $a_j\not=0$ for some $j$. Put $N=\min\{j\geq 1: a_j\not=0\}$ and $q=\fr{a_N}{b_{N}}$.  Note $q\in\T$. By (\ref{e2}), we have
$
\fr{\al_{m+N}}{\al_m}=q
$
for all  $m\geq 0$, which yields
\begin{equation}\label{e3}
\fr{\al_{m+Nk}}{\al_m}=q^k, 
\qquad m\geq 0, \ k\ge 1.
\end{equation}
Now, in order to see that (c) holds, we consider the following two cases.

{\it Case 1.} Suppose $a_{Nk}\not=0$ for some $k\ge 1$. By (\ref{e2}) and (\ref{e3}), we have
$$
\fr{a_{Nk}}{b_{Nk}}=\fr{\al_{Nk}}{\al_0}=q^k
$$
for all $k\ge  1$ with $a_{Nk}\not=0$. Hence $b_{Nk}=\bar q^ka_{Nk}$ for all $k\ge 1$.

{\it Case 2.}
Suppose 
$a_{Nk+j}\not=0$ for some $k\ge 1$ and $j$ with $1\le j<N$.
By (\ref{e2}) and (\ref{e3}) again, we also see that
$$
\fr{a_{Nk+j}}{b_{Nk+j}}=\fr{\al_{Nk+j+\ell}}{\al_\ell}=q^{k}\fr{\al_{j+\ell}}{\al_\ell}
$$
for $0\le\ell\le N-j-1$ and
$$
\fr{a_{Nk+j}}{b_{Nk+j}}
=\fr{\al_{N(k+1)+\ell-(N-j)}}{\al_\ell}
=q^{k+1}\fr{\al_{\ell-(N-j)}}{\al_\ell}
$$
for $N-j\le\ell< N$. It follows that
\begin{align*}
\fr{\al_{j}}{\al_0}
&=\fr{\al_{j+1}}{\al_1}
=\cdots
=\fr{\al_{N-1}}{\al_{N-j-1}}\\
&=q\fr{\al_0}{\al_{N-j}}
=q\fr{\al_1}{\al_{N-j+1}}
=\cdots
=q\fr{\al_{j-1}}{\al_{N-1}}.
\end{align*}
Multiplying all terms of the above together, we obtain
$$
\fr{\al_{j}}{\al_0}=q^{j/N},
$$
which gives
$$
\fr{a_{Nk+j}}{b_{Nk+j}}=q^{k}\fr{\al_{j}}{\al_0}=q^{(Nk+j)/N}.
$$
Therefore, the above shows that $b_{Nk+j}=\bar q^{(Nk+j)/N}a_{Nk+j}$ holds for any $a_{Nk+j}\ne 0$, so does for any $a_{Nk+j}$ with $k\ge 1$ and $1\le j<N$.

Noting $q\in \T$ and letting $\lam=\bar q^{1/N}\in\T$, we see from the two cases above, $b_j=\lambda^{j}a_j$ holds
for all integers $j\ge 1$.
Therefore we conclude that (c) holds with 
$$
f=\fr{a_0+b_0}{2}+\sum_{j=1}^\infty a_jz^j\in H^2.
$$
 The proof is complete.
\end{proof}

\begin{rk}
We remark that it may not hold $\al_j=q^{j/N}\al_0$ for any $j$; see
Example \ref{e-5-23}. We also note that the argument above proving $b_j=\lambda^j a_j$ is different from and much simpler than that in \cite[Proposition 2.1]{BCZ} where canonical symmetric Toeplitz operators have been characterized.
\end{rk}


\begin{rk}\label{rk-1}
Let $\phi\in L^\infty(\T)$, $\lam\in\T$ and $u$ be an inner function with $u(a)=0$ for some $a\in\D$.
Choose an \onb~ $E=\{e_a,,e_1,e_2,\ldots\}$ of $H^2\ominus uH^2$. If $T_\phi$ is $\Clae$-symmetric,
by the similar argument as in the proof of Theorem \ref{C1-1}, we see that there are $f, g_j\in H^2 (j\geq 1)$ such that
$$
\phi=f\circ u+\ov{C_\lam(f)\circ u} +\sum_{j\geq 1} \left[\frac{e_j}{e_a}\cdot g_j\circ u+\ov{\frac{e_j}{e_a}\cdot C_\lam(g_j)\circ u} \, \right].
$$
\end{rk}

Now we present examples asserting that 
$C_{\lam,u,E}$-symmetry of a Toeplitz operator generally depends on a choice of \onb \ $E$.
In addition, the example also exhibits that a Toeplitz operator with the symbol $\phi$ appeared in the preceding remark
may be complex symmetric.

\begin{ex}\label{p-7}
Given $a\in\D\setminus\{0\}$, we let
$u=z\va_a$ and consider $E=\{1,e\}$ where  $e=\gamma\fr{K_a-1}{||K_a-1||}$ and $\gamma\in\T$. Let $\lam=\fr{a^2}{\bar a^2}$ and $\phi\in L^\infty(\T)$. Using Theorem \ref{C1-1} and Remark \ref{rk-1},
one can check that $T_\phi$ is $\Clae$-symmetric if and only if exactly one of the following statements holds:
\begin{enumerate}
\item[(a)]  $\gamma\not=\pm 1$ and $\phi=f\circ u+\ov{C_\lam(f)\circ u}$ for some $f\in H^2$.
\item[(b)] $\gamma=\pm 1$ and
$\phi=f\circ u+\ov{C_\lam(f)\circ u}+e(g\circ u)+\ov{e[ C_\lam(g)\circ u]}$ for some $f,g\in H^2$.
\end{enumerate}
\end{ex}


In the rest of this section,  let's consider an inner function given by a single M\"obius transform.
To be more precise, given $a\in \D$, $\lambda\in\T$ and $\alpha\in\mathscr M$, let $e_a=\sqrt{1-|a|^2}K_a$.
Consider $u=\va_a$ and $E_\gamma=\{\gamma e_a\}$ where  $\gamma\in\T$.
We put
$$
C_{\alpha,a}:=C_{\alpha,\va_a,E_1}, \qquad C_{\lam,a}:=C_{\lam,\va_a,E_1}
$$
for notational simplicity.
Since $C_{\alpha,\va_a,E_\gamma}=\gamma^2 C_{\alpha,a}$, as an immediate consequence of  Proposition
\ref{C1} and Theorem \ref{C1-1}, we have the following corollary which recovers known results in \cite{BCZ} and \cite{KL} in case $a=0$.

\begin{cor}\label{cor-1}
Let $\phi\in L^\infty(\T)$ and $a\in\D$. Then the following are equivalent:
\begin{itemize}
\item[(a)] $T_\phi$ is $C_{\alpha,a}$-symmetric for some $\alpha\in\mathscr M$.
\item[(b)] $T_\phi$ is $C_{\lam,a}$-symmetric for some $\lam\in\T$.
\item[(c)] $T_{\phi\circ\va_a}$ is $C_\lam$-symmetric for some $\lam\in\T$.
\item[(d)] $\phi=f\circ\va_a+\ov{C_\lam(f)\circ\va_a}$ for some $\lam\in\T$ and $f\in H^2$.
\item[(e)] $\phi\circ\va_a(\lam z)=\phi\circ\va_a(\bar z)$ for all $z\in \D$ and some $\lam\in\T$.
\end{itemize}
\end{cor}

%


Given $a\in \D$ and $\lambda\in \T$, since 
\begin{equation*}
\label{e-cla}
U_aC_\lam U_a(\va_a^me_a)=\bar\lam^m\va_a^me_a
=C_{\lam,a}(\va_a^me_a)
\end{equation*}
for all integer $m\ge 0$,
we have $C_{\lam,a}=U_aC_\lam U_a$ 
Using this fact, for $\phi\in L^\infty(\T)$, we see that  $T_\phi$ is $C_{\lam,a}$-symmetric if and only if
$U_aT_\phi U_a$ is $C_{\lam}$-symmetric. Noting $U_aT_\phi U_a=T_{\phi\circ\varphi_a}$, we  provide another proof of the equivalences of (b), (c), (d) and (e) of Corollary \ref{cor-1}.

We give in passing an example of a Toeplitz operator which is complex symmetric with respect to 
 $C_{\alpha,0}$ and $C_\lambda$ for all $\lambda\in \T$ and some $\alpha\in \mathscr M$, although $C_{\alpha,0}\ne \mu C_{\gamma,a}$ for any  $a\in\D\setminus\{0\}$
and $\mu,\gamma\in\T$.

\begin{ex}\label{e-5-23}
Let $\lam\in\T$ and  $\eta\in\T\setminus\{\pm\bar\lam\}$. Define a conjugation $C$ on $H^2$ by
$$
C z^{2j}=\bar\lam^{2j}z^{2j}, \qquad C z^{2j+1}=\eta\bar\lam^{2j}z^{2j+1}
$$
for $j=0,1,2,\cdots.$ So $C=C_{\alpha,0}$ where
$$
\alpha=(1,\eta,\bar\lam^2,\eta\bar\lam^2,\bar\lam^4,\eta\bar\lam^4,\cdots)\in\mathscr M.
$$
Note $C\ne \mu
C_\gamma$. We claim $C\ne \mu C_{\gamma,a}$ for any  $a\in\D\setminus\{0\}$
and $\mu,\gamma\in\T$. Suppose not. Then, $C=\mu C_{\gamma,a}$ for some $a,\mu,\gamma$. Then, a direct calculation shows 
$$
C z^n=\mu \delta (\beta\va_b)^ne_b,\qquad n=0,1,2,\cdots
$$
where $b=\va_{a}(\bar a\gamma)$, $\delta=\fr{1-a^2\bar\gamma}{|1-a^2\bar\gamma|}$ and $\beta\in\T$ satisfying
$\bar b=\beta b$. Comparing the above with the definition of $C$, we have $b=0$ and then $\gamma=a/\bar a$.
Hence $C=\mu C_{-\ov\beta}$ and  a contradiction follows.
Then, for $\phi=z^2+z^4+\lam^2\bar z^2+\lam^4\bar z^4$, it is clear that $T_\phi$ is complex symmetric with respect to both $C_\lam$ and $C$.
\end{ex}

We now provide another point of view for conjugations $C_{\lambda,a}$.
Let $u: \D\to\C$ and $v: \D\to\D$ be two analytic functions. Define a conjugate-linear operator $\mathcal{A}_{u,v}$ on $H^2$ by
$$
\mathcal{A}_{u,v}f=u\ov{f(\ov{v})}
$$ 
for functions $f\in H^2$.
Such conjugations are introduced in \cite{LimK} where complex symmetric weighted composition operators have been studied.
Clearly, if $u=\mu$ and $v=\bar\lam z$ for some $\mu, \lam\in\T$, then $\mathcal{A}_{u,v}=\mu C_\lam$.
So $\mathcal{A}_{u,v}$ is a generalization of the conjugations $\mu C_\lam$.

We prove that an operator
$\mathcal{A}_{u,v}$ being a conjugation is essentially  the same as the conjugations $C_{\lam,a}$ as shown in the following proposition, which is not only of own interest but also used 
in Proposition \ref{p-new-1}.

\begin{prop}\label{p-add}
Let $u: \D\to\C$ and $v: \D\to\D$ be two analytic functions. Then the following  assertions are equivalent:
\begin{enumerate}
\item[(a)] $\mathcal{A}_{u,v}$ is a conjugation on $H^2$.
\item[(b)] $u=\alpha e_b$ and $v=\beta\va_b$ for some $\alpha,\beta\in\T$ and $b\in\D$ with $\bar b=\beta b$.
\item[(c)] $\mathcal{A}_{u,v}=\mu C_{\lam,a}$ for some $\mu,\lam\in\T$ and $a\in\D$.
\end{enumerate}
\end{prop}

\begin{proof}
The equivalence of (a) and (b) has been proved in \cite{LimK}.
Now assume (c) and show (b). 	
If  $a=0$, then (b) holds with $u=\mu$ and $v=\bar\lam z$.
So, suppose $a\not=0$. If  $\lam\not=a/\bar a$, then a  straightforward calculation yields
$\ov{e_a(\lam\ov{\va_a})} e_a=\ov{\gamma}e_{\va_a(\bar a\lam)}$ and 
\begin{align}
\label{e-13000}
\ov{\va_a(\lam\ov{\va_a})}
=\va_{\ov{a^2}}(\bar\lam)\va_{\va_a(\bar a\lam)}
=(\bar b/b)\va_b
\end{align}
where
\begin{equation*}
b=\va_a(\bar a\lam),
\qquad 
\gamma=\fr{1-a^2\bar\lam}
{|1-a^2\bar\lam|}
	\end{equation*}
	for simplicity.
It follows that
\begin{align}
\label{e-130}
\begin{aligned}
C_{\lam,a}z^n
&=U_aC_\lam U_az^n\\
&=\ov{[\va_a(\lam\ov{\va_a})]^n}\ov{e_a(\lam\ov{\va_a})} e_a \\
&=\ov{\gamma}[(\bar b/b)\va_b]^n e_b
\end{aligned}
\end{align}
for all integer $n\ge 0$.
	Therefore, we have
	$$
	C_{\lam,a}f=\ov{\gamma} e_b \ov{f(\ov{(\bar b/b)\va_b})}
	$$
for all $f\in H^2$.	So (b) holds with $u=\mu\ov{\gamma} e_b$ and $v=(\bar b/b)\va_b$.
If $\lam=a/\bar a$, it follows from (\ref{e-130}) that
	$
	C_{\lam,a}z^n=(\bar\lam z)^n 
	$ for all integer $n\ge 0$. Thus
	$
	C_{\lam,a}f=\ov{f(\ov{\bar\lam z})}
	$ for all $f\in H^2$ and (b) follows with $u=\mu$ and $v=\bar\lam z$.

Finally assume (b) and show (c).	
	If $b=0$, then it is obvious that $\mathcal{A}_{u,v}=\alpha C_{-\ov\beta}$.
	Now suppose $b\not=0$. It is routine to see that there exist $\lam\in\T$ and $a\in\D$ such that $b=\va_a(\bar a\lam)$.
Put $\gamma:=\fr{1-a^2\bar\lam}{|1-a^2\bar\lam|}$. Then,  (\ref{e-130}) shows that
	$$
	C_{\lam,a} z^n=\ov{\gamma}(\beta\va_b)^n e_b,\qquad n=0,1,2,\cdots
	$$
because $\beta=\bar b/b=\va_{\ov{a^2}}(\bar\lam)$. Thus $\mathcal{A}_{u,v}=\alpha\gamma C_{\lam,a}$ and (c) holds. 
The proof is complete.
\end{proof}


\section{A class of \uet \ Toeplitz operators}\label{S:Blaschke}

In this section, we first show that Question \ref{C:main} is true for certain Toeplitz operators induced by a M\"obius transform. We also
prove Theorem \ref{T:main-2}.

A conjugate-linear map $D$ on a Hilbert space $\HH$ is called  anti-unitary if $D$ is invertible and isometric.
Clearly, a conjugation is anti-unitary.
Recall that an operator $T$ on  $\HH$ is called \uet\ if $T\cong CT^*C$ for some conjugation $C$ on $\HH$.
The following lemma shows that \uet~ operators  can be characterized by anti-unitary operators; see \cite{GJZ} for details.

\begin{lem}
For a bounded linear operator $T$ on $\HH$,
$T$ is \uet \ if and only if $DTD^{-1}=T^*$ for some anti-unitary operator $D$ on $\HH$.
\end{lem}


We introduce nowhere
winding curves which are closely related to complex symmetric/\uet\ Toeplitz operators.
For a closed curve $\gamma$ in $\mathbb C$,   
the winding number of $\gamma$ with respect to $z\in\C\setminus \gamma$, also called the  index of $\gamma$,
is defined by
$$
{\rm Ind}_{\gamma}(z)=\fr 1{2\pi i}\int_{\gamma} \fr{1}{\zeta-z}\,d\zeta,\qquad z\in\C\setminus \gamma.
$$
Recall that the curve $\gamma$ is  called nowhere winding if ${\rm Ind}_{\gamma}(z)=0$ for all
$z\in\mathbb C\setminus \gamma$.
Given $\phi\in L^\infty(\T)$, a well known Coburn theorem  says that at least one of $T_\phi$ and $T^*_\phi$ is injective; see \cite{Co} or \cite{Do} for details. Using this fact together with the Fredholm index theorem(\cite{Do}), Noor \cite{Noor} recently  proved that
the curve $\phi(\T)$ of the continuous symbol $\phi$ of a complex symmetric Toeplitz operator $T_\phi$ is nowhere winding.
Note each complex symmetric operator is an \uet~ operator.
By using a similar argument as in  \cite{Noor}, we can see that the same is true for \uet~ operators as shown in the following.

\begin{prop}\label{p-4}
Let $\phi$ be a continuous function on $\T$. If $T_\phi$ is \uet, then
$\phi(\T)$ is a nowhere winding curve.
\end{prop}

In \cite{Noor}, Noor posed a question: {\it Does every nowhere winding curve 
induce a complex symmetric operator?} In view of Proposition \ref{p-4}, it is natural to ask the same question for \uet \ operators.
 The following provides a positive answer 
for certain Toeplitz operators.

\begin{prop}
\label{th-4}
For $q\in\T$ and  $b\in\D$, let $\phi=z+\ov{q\va_b}$.
 Then the following statements  are equivalent:
\begin{enumerate}
\item[(a)] $T_\phi$ is complex symmetric.
\item[(b)] $T_\phi$ is \uet.
\item[(c)] $T_\phi$ is either normal or $C_{\lam,a}$-symmetric for some $\lam\in\T$ and $a\in\D$.
\item[(d)] The curve $\phi(\T)$ is nowhere winding.
\item[(e)] $qb=\bar b$.
\end{enumerate}
\end{prop}

\begin{proof}
Since implications (c) $\Longrightarrow$ (a) $\Longrightarrow$ (b) is obvious and 
(b) $\Longrightarrow$ (d) follows from
Proposition \ref{p-4}, we only need to prove implications (d) $\Longrightarrow$ (e) and (e) $\Longrightarrow$ (c).
First assume (d) and show (e).
Put $s:=\bar q\bar b-b$. By a direct computation, we see
$$
\phi(z)-s=\fr{z^2-[\bar q(1-|b|^2)+b^2]}{z-b}, \qquad z\in\T.
$$
If $qb\not=\bar b$, a routine calculation shows
$
|\bar q(1-|b|^2)+b^2|<1,
$
which implies that $s\notin \phi(\T)$ and two zeros of  $z^2-[\bar q(1-|b|^2)+b^2]$ are all in $\D$. Hence, by the argument principle, we see that
$$
{\rm Ind}_{\phi(\T)}(s)=\fr 1{2\pi i}\int_{\phi(\T)}\fr{1}{\zeta-s}\, d\zeta
=\fr 1{2\pi i}\int_\T\fr{\phi'(\zeta)}{\phi(\zeta)-s}\, d\zeta=1,
$$
which is a contradiction and  (d) $\Longrightarrow$ (e) holds.

Finally, assume (e) and show (c).  
If $b=0$, then $\phi=z-\ov{qz}$ and $T_\phi$ is
a normal operator.
Now suppose $b\not=0$. It is routine to see that there exist $a\in\D$ and $\lam\in\T$ such that
$b=\va_a(\bar a\lam)$. Then, using  applications of (\ref{e-new-ee2}) and (\ref{e-13000}), we can check that
$$
C_{\lam,a}T_{z+\ov{q\va_b}}C_{\lam,a}=T_{C_\lam(\va_a+\ov{q\va_b\circ\va_a})\circ\va_a}
=T^*_{z+\ov{q\va_b}},
$$
thus $T_{z+\ov{q\va_b}}$ is $C_{\lam,a}$-symmetric and (c) follows, as desired.
\end{proof}


Given an operator $T$, we let $\ran\,T$, $\ker T$ and ${\rm tr}\, T$ denote the range, kernel and trace of $T$ respectively.
Also, let $[A, B]=AB-BA$ denote the commutator of two operators $A$ and $B$.

The following lemma will be useful in the proof of Theorem \ref{T:main-2}.
Given an inner function $u$, we denote ${\rm ord}\, u=\dim H^2\ominus uH^2$.

\begin{lem}\label{p-2}
Let $f, g$ be finite Blaschke products and  $\mu\in\C$. Put $\phi=f+\mu\bar g$. If $DT_\phi=T^*_\phi D$
for some anti-unitary operator $D$ on $H^2$, then the following statements hold:
\begin{enumerate}
\item[(a)] $[T_\phi^*, T_\phi]$ has finite rank and ${\rm ran}\,[T_\phi^*, T_\phi]$ reduces $D$.
\item[(b)]  ${\rm ord}~f=|\mu|^2 {\rm ord}~g.$
\end{enumerate}
\end{lem}

\begin{proof}
Set $T=[T^*_\phi, T_\phi]$.
It is easy to check
\begin{equation}\label{e-commu}
T=[T^*_f, T_f]-|\mu|^2 [T^*_g, T_g].
\end{equation}
Since $T_f$ is an isometry, $[T^*_f, T_f]=I-T_fT^*_f$ is a projection onto $H^2\ominus fH^2$ and hence
$$
{\rm tr}[T^*_f, T_f]=\dim (H^2\ominus fH^2)={\rm ord}\,f.
$$
Similarly, we have
$
{\rm tr}[T^*_g, T_g]={\rm ord}\, g.
$
Thus $T$ is a finite-rank operator and ${\rm tr}\,T={\rm ord}~f-|\mu|^2 {\rm ord}~g$. Also,
since $DT_\phi =T_\phi^*D$, it follows that $DT=-TD$. Thus $D(\ker T)\subset\ker T$ and $D(\ran\, T)\subset\ran\, T$.
Noting that $T$ is self-adjoint, we see that $\ker T\oplus \ran\, T=H^2$ and hence $\ran\, T$ reduces $D$. This gives (a).

The condition $DT_\phi=T^*_\phi D$ also yields $DTD^{-1}=-T$. It follows that ${\rm tr}\,T=0$, which implies (b). The proof is complete.
\end{proof}

Now we are ready to prove Theorem \ref{T:main-2}.

\begin{proof}[Proof of Theorem \ref{T:main-2}]
Set $\psi:=z+\ov{\mu B_2}$ and then $\phi=\psi\circ B_1$.
Let ${\rm ord}\,B_1=m$.

First assume (b) and show (c). Since $T_\phi\cong T_\psi^{(m)}$ (see \cite{Cowen}) and $T_\psi$ is essentially normal, by the similar argument as in the proof of
Theorem 1.4 of \cite{BCZ}, we see $T_\psi$ is also \uet. By Proposition \ref{th-4},
it suffices to prove that $\mu\in\T$ and ${\rm ord}\, B_2=1$.
Suppose $D^{-1}T_\psi D=T^*_\psi$ for some anti-unitary $D$ on $H^2$. By (\ref{e-commu}) and Lemma \ref{p-2}(b), we have
$$
T:=[T^*_\psi, T_\psi]=(I-T_zT_{\bar z})-|\mu|^2(I-T_{B_2}T_{\bar{B_2}})
$$
and 
$
1=|\mu|^2 {\rm ord}\, B_2,
$
which means that $0<|\mu|\leq 1$. If $|\mu|=1$, then ${\rm ord}\, B_2=1$, as desired.
So, suppose $|\mu|<1$ and put $\mathscr R={\rm ran}\,T$.
Noting that $I-T_zT_{\bar z}$ is a projection onto $\C$ and $I-T_{B_2}T_{\bar B_2}$ is a projection onto $H^2\ominus B_2H^2$, we have
\begin{equation}
\label{e-subset}
\mathscr R\subset\vee\{1, H^2\ominus B_2 H^2\}.
\end{equation}
We claim $1\in\mathscr R$. In fact, if $B_2(0)=0$, then $T1=1-|\mu|^2$; if $B_2(0)\not=0$, then
$TB_2=B_2(0)$, which proves the claim.
Now, by Lemma \ref{p-2}(a), we know that $D(\mathscr R)\subset\mathscr R$. Since $1\in\mathscr R$, in view of (\ref{e-subset}),
we may assume
$
D1=\alpha+h
$
for some $\alpha\in\C$ and $h\in H^2\ominus B_2H^2$. Recalling 
$
T^*_\psi 1=D^{-1}T_\psi D1
$
and noting $H^2\ominus B_2H^2=\ker T_{\bar B_2}$, we obtain  
$$
\mu B_2=T^*_\psi 1=D^{-1}[zD1+\bar\mu\alpha\ov{B_2(0)}\,]
$$
which implies
$$
|\mu|^2=||\mu B_2||^2=||zD1+\bar\mu\alpha\ov{B_2(0)}||^2=1+|\mu\alpha B_2(0)|^2.
$$
So we have
$$
1=|\mu|^2(1-|\alpha|^2|B_2(0)|^2),
$$
which contradicts the hypothesis $|\mu|<1$. Then we conclude that $|\mu|=1$ and ${\rm ord}\, B_2=1$ as desired.

Now assume (c).  By Proposition \ref{th-4}, $T_\psi$ is complex symmetric, then so is $T_\phi$ because
$T_\phi\cong T_\psi^{(m)}$. Hence (a) holds.
Since implication (a)$\Longrightarrow$(b) is obvious, we complete the proof.

\end{proof}

\section{Analytic symmetry of Toeplitz operators}\label{S:Regular}

In this section, we first show that
 an analytic symmetric Toeplitz operator with trigonometric symbol is $C_\lambda$-symmetric for some $\lambda \in\T$.
 Also, we characterize
certain Toeplitz operators to be analytically symmetric.

Given $\phi\in L^\infty(\T)$,  we recall that $T_\phi$ is  analytically symmetric if there exist $f,g\in H^\infty$ and a conjugation $C$ on $H^2$ such that $CT_fC=T_g$ and $\phi=f+\bar g$.
Also, analytic symmetric Toeplitz operators are complex symmetric.


To begin with, for $f\in H^\infty$, inner $u$ and $\lam\in\T$,
we first recall $C_\lambda f(z)=\overline{f(\lambda \bar z)}$  and note from (\ref{e-new-ee2}) that the Toeplitz operator of form
$$
T_{f\circ u+\ov{C_\lam (f)\circ u}}=T_{f\circ u+\ov{C_1(f)\circ(\bar\lam u)}}
$$ is analytically symmetric.
The following lemma shows that 
any analytic symmetric Toeplitz operator
has the similar form above under certain condition.

\begin{lem}\label{lem-new}
Let $f,g\in H^\infty$ and $\alpha\in\C$.
Suppose the inner factor
of $f-\alpha$ is a nonconstant finite Blaschke product. If $CT_fC=T_g$ for some conjugation
$C$, then there exist $h\in H^\infty$
and finite Blaschke products $u, u_1$ with ${\rm ord}\, u={\rm ord}\, u_1$ such that
$
f=h\circ u$ and  $g=C_1h\circ u_1.
$
\end{lem}

\begin{proof}
Put $U=CC_1$. Then, it is easy to see that $U$ is a unitary operator and 
$$U^*T_fU=C_1CT_fCC_1=C_1T_gC_1=T_{C_1g}$$ on $H^2$. It follows from \cite[Corollary 1]{Cowen} that
	there exist $h\in H^\infty$ and two finite Blaschke products $u, u_2$ with ${\rm ord}\, u={\rm ord}\, u_2$ such that
	$
	f=h\circ u$ and $ C_1g=h\circ u_2$,
	which immediately  implies the desired result  by taking $u_1=C_1u_2$. We complete the proof.
\end{proof}

For trigonometric symbols $\phi$, we show below that Toeplitz operator $T_\phi$ is analytically symmetric
only when $\phi=f+\ov{C_\lam f}$ where $f$ is an analytic polynomial.

\begin{thm}\label{p-6}
	Let $\phi$ be a trigonometric polynomial. Then $T_\phi$ is analytically symmetric
	if and only if $T_{\phi}$ is $C_\lam$-symmetric for some $\lam\in\T$.
\end{thm}

\begin{proof}
By Proposition \ref{C1}, the sufficiency is clear because $\phi=f+\ov{C_\lam f}$ and $C_\lam T_f C_\lam=T_{C_\lam f}$
for all $\lambda\in\T$ and $f\in H^\infty$.

Now we prove the necessity.	
	Suppose $\phi=f+\bar g$ and $CT_fC=T_g$ for some conjugation $C$ and
	 analytic polynomials $f,g$. Since $f-\alpha$ is still an analytic polynomial for any $\alpha\in\C$, 
 its inner factor is a finite Blaschke product. Thus by Lemma \ref{lem-new},
	there exist $h\in H^\infty$ and two finite Blaschke products $u$ and $u_1$ with $ {\rm ord}\, u={\rm ord} \, u_1$ such that
	$$
	f=h\circ u, \ \ \ g=(C_1h)\circ u_1.
	$$
	Since $f$ and $g$ are both analytic polynomials, we may assume that $h$ is an analytic polynomial and
	$u=\beta z^N$, $u_1=\beta_1 z^N$ where $\beta,\beta_1\in\T$ and $N$ is a positive integer.
	Hence
	$$
	f+\bar g=h\circ(\beta z^N)+\ov{C_1(h)\circ(\beta_1 z^N)}=h\circ  u+\ov{C_\lam (h)\circ  u},
	$$
	where $\lam=\beta\ov{\beta_1}$. Thus by Proposition \ref{C1}, $T_{f+\bar g}$
	is $\Clae$-symmetric where $E=\{1, z, z^2,\cdots, z^{N-1}\}$. It follows from the definition of $\Clae$  that
	$$
	\Clae( u^m z^j)=\bar\lam^m u^m z^j, \qquad m\geq 0, \ j=0,1,\cdots,N-1.
	$$
	Then,  by Corollary \ref{cor-1}, $T_{f+\bar g}$ is $C_\lam$-symmetric for some $\lam\in\T$, as
	desired. The proof is complete.
\end{proof}

We provide an example of a complex symmetric Toeplitz operator which is not analytically symmetry.

\begin{ex}\label{E:counterEg}
	Let $\phi=z+z^2+i(\overline{ z+z^2})$. Then $T_\phi$ is obviously  normal and hence complex symmetric.
	However, by Proposition \ref{C1}(f), we see that $T_\phi$ is not $C_\lam$-symmetric
	for any $\lam\in\T$ and hence $T_\phi$ is not analytically symmetric  by Theorem \ref{p-6}.
\end{ex}

In the following theorem, we characterize
an analytic symmetric Toeplitz operator when the analytic part of its symbol is a polynomial of degree 2.
We were not able to characterize general symbols.

\begin{thm}
\label{p-new-1}
Let $\ell, c, d\in\mathbb C$ and $g\in H^\infty$ with $g(0)=0$. Put $p=dz+cz^2$ and
$\phi=\ell+p+\bar g$. 
Then the following statement are equivalent:
\begin{enumerate}
\item[(a)] $T_\phi$ is analytically symmetric on $H^2$.
\item[(b)] There exist $q\in\T$ and $b\in \D$ with $\bar b=qb$ such that one of the following statements holds:
\begin{enumerate}
\item[(b1)] $d\ne 0$ and $g=\bar dq\va_{b}+\bar c(q\va_{b})^2$.
\item[(b2)] 
$d=0$ and  $g=\bar c(q\va_b)^2$.
\item[(b3)] 
$d=0$ and $g=\bar c q\va_{b}(z^2)$.

\end{enumerate}
\end{enumerate}
\end{thm}

\begin{proof} 
We first assume (b) and prove (a).
By (\ref{e-13000}), we note 
$C_\lam(\va_a)\circ\va_a=q\va_b$ 
where $\bar b=qb$ with $b=\va_a(\bar a\lam)$. It follows from
(\ref{e-new-ee2}) that	
$$
C_{\lam,a}T_pC_{\lam,a}
=T_{C_\lam(p\circ\va_a)\circ\va_a}=T_{\ov{p(\bar q\ov{\va_b})}}
$$
	and 
\begin{equation*}
\label{e-new-1-3}
		C_{\lam, \va_a\circ z^2, E}T_{z^2} C_{\lam, \va_a\circ z^2, E}=T_{C_\lam(\va_a)\circ\va_a\circ z^2}=T_{q\va_b\circ z^2}
	\end{equation*}
for an \onb \ $E$ of the model space $H^2\ominus(\va_a\circ z^2)H^2$.
	Hence Toeplitz operators $T_{p+p(\bar q\ov{\va_b})}$ and $T_{z^2+\ov{q\va_b\circ z^2}}$ with $\bar b=qb$ both are
	analytically symmetric and (a) holds.
	
	Now we assume (a) and prove (b). We suppose that
	\begin{equation}\label{new-e-2}
		CT_{p+\ell_1}C=T_{g+\bar\ell_2}
	\end{equation}
	for some conjugation $C$, here $\ell_1+\ell_2=\ell$. Since $p+\ell_1$ is an analytic polynomial, by Lemma \ref{lem-new},
	one may assume
	$$
	g+\bar\ell_2=C_1(p+\ell_1)\circ B=\bar\ell_1+C_1(p)\circ B
	$$
	for some finite Blaschke product $B$ with ${\rm ord}B=1$  or $2$ according to $\ell_1, d$ and $c$.  Thus (\ref{new-e-2}) becomes
	$
	CT_{p}C=T_{g}$ and $g=C_1(p)\circ B.
	$
	Hence we easily see
	\begin{equation}\label{e-new-1}
		C(\ker T^*_{p-\beta})=\ker T^*_{g-\bar\beta}
	\end{equation}
for any $\beta\in\C$. 

Now we consider several cases.  First suppose $d\not=0$ and $c\neq 0$. In this case, we may take $B=q\va_b$ for some $q\in\T$ and $b\in\D$ and hence
	$
	g=\bar dq\va_{b}+\bar c(q\va_{b})^2.
	$
Put
	$$
	E_d=\{z\in\D: \ |(d/c)+z|>1\}.
	$$
	Then $E_d$ is a nonempty  open set  in $\D$ because $d/c\ne 0$. 
Let $z\in E_d$ be any point. Then we note
	$$
	\ker T^*_{p-\beta}=\C\cdot K_z, \qquad \ker T^*_{g-\bar\beta}=\C\cdot K_{\va_{b}(\bar q\bar z)}
$$
where $\beta:=dz+cz^2$.
	Hence, by (\ref{e-new-1}), one can see that there is $\eta(z)\in \C$ such that
	$
	CK_z=\eta(z) K_{\va_{b}(\bar q\bar z)}.
	$
	It follows that  
	\begin{align*}
		Ch(z) &=\langle Ch, K_z\rangle=\langle CK_z, h\rangle=\eta(z)\ov{h(\ov{ q\va_{\bar q\bar b}(z)})}
	\end{align*}
	for every $h\in H^2$. By considering $h=1$,  we have $\eta=C1$ and hence
	$$
	Ch(z)=C1(z)\ov{h(\ov{q\va_{\bar q\bar b}(z)})},\qquad h\in H^2
	$$
	for each $z\in E_d$. Because $E_d$ is open in $\D$, the above holds for any $z\in\D$. Therefore, an application of 
	Proposition \ref{p-add} shows $\bar b=qb$ and then (a) holds.

If $d\ne 0$ and $c=0$, one can see that (a) holds by the similar proof.

Now, suppose $d=0$ and assume $p=z^2$ without loss of generality. As before we may
further	assume $g=(q\va_b)^2$ or $g=q\va_a\va_b$, where $q\in\T$ and $a,b\in\D$
	with $a\ne b$. We divide the proof into two cases below.
	
{\it Case 1}. Consider $g=(q\va_b)^2$. Since the case $b=0$ is trivial, we assume $b\ne 0$.
	Taking $\beta=0$ in (\ref{e-new-1}) we have
	\begin{equation}\label{e-new-4}
		C1=u_1K_b+v_1\partial_{\bar b}K_b
	\end{equation}
	for some $u_1,v_1\in\C$ where $\partial_{\bar b}K_b(z)=z(1-\bar b z)^{-2}$. Again, taking $\bar\beta=q^2b^2$ in
	(\ref{e-new-1}) we see $1\in\ker T^*_{g-\bar\beta}$ and $K_{\bar q\bar b}, K_{-\bar q\bar b}\in\ker T^*_{p-\beta}$,
	so there are $u_2,v_2\in\C$ such that
	$$
	C1=u_2K_{\bar q\bar b}+v_2K_{-\bar q\bar b}.
	$$
	Comparing the above with (\ref{e-new-4}), we obtain $b=\bar q\bar b$ or $b=-\bar q\bar b$. Hence (b) follows, as desired.
	
{\it Case 2}. Consider $g=q\va_a\va_b$ with $a\ne b$. We first assume $ab\ne 0$. Taking $\beta=0$ in the equality (\ref{e-new-1}) we have
	\begin{equation}\label{e-new-2}
		C1=u_1K_a+v_1K_b, \ \ \ Cz=u_2K_a+v_2K_b
	\end{equation}
	for some $u_j,v_j\in\C, j=1,2$. Taking $\bar\beta=qab=:\bar\al^2$, then $1\in \ker T^*_{g-\bar\beta}$ and $K_\al, K_{-\al}
	\in\ker T^*_{p-\beta}$. So there are $u_3,v_3\in\C$ such that
	$$
	C1=u_3K_\al+v_3K_{-\al}.
	$$
	Comparing the above with $C1$ in (\ref{e-new-2}), we may assume $\al=a$ without loss of generality, and hence $\bar a^2=qab$
	and $u_1=u_3$. Now we will  show $b=-a$. Then we have
	$g=q\va_a\va_{-a}=-q\va_{a^2}\circ z^2$ and (b) holds.
	If $v_1\not=0$ or $v_3\not=0$, then we get $b=-\al=-a$, as desired. Now we suppose $v_1=v_3=0$. Hence (\ref{e-new-2}) becomes
	\begin{equation*}
		C1=u_1K_a, \qquad Cz=u_2K_a+v_2K_b.
	\end{equation*}
Since $CT_{z^2}(z^n)=T_gC(z^n)$ for all integer $n\ge 0$, it follows that
	\begin{equation}\label{e-new-3}
		Cz^{2n}=g^n u_1 K_a, \qquad Cz^{2n+1}=g^n(u_2K_a+v_2K_b)
	\end{equation}
for all integer $n\ge 0$.
Put $W=CC_1$. Then $W$ is unitary and $C_1WC_1=W^*$. Also,
	$$
	Wz^{2n}=g^nu_1 K_a, \ \ \ Wz^{2n+1}=g^n(u_2K_a+v_2K_b)
	$$
for all integer $n\ge 0$.
	So, for each $h\in H^2$ with $h(z)=\sum_{n=0}^\infty a_nz^n$, we have
	\begin{align*}
		Wh(z) &=\sum_{n=0}^\infty a_{2n} W(z^{2n})+\sum_{n=0}^\infty a_{2n+1}W(z^{2n+1})\\
		&=u_1K_a\sum_{n=0}^\infty a_{2n}g^n+(u_2K_a+v_2K_b)\sum_{n=0}^\infty a_{2n+1}g^n\\
		&=U_1(z) h(V_1(z))+U_2(z) h(V_2(z))
	\end{align*}
	for $z$ in a connected open set $E\subset\D-\{a,b\}$ where
	$$
	U_1=\fr 1{2}\Big(u_1K_a+\fr{u_2K_a+v_2K_b}{\sqrt{g}}\Big), \qquad V_1=\sqrt{g},
	$$
	$$
	U_2=\fr 1{2}\Big(u_1K_a-\fr{u_2K_a+v_2K_b}{\sqrt{g}}\Big), \qquad V_2=-\sqrt{g}.
	$$
Note
	$$
	\langle W^*K_z,h\rangle=\langle K_z, Wh\rangle=\langle \ov{U_1(z)}K_{V_1(z)}+\ov{U_2(z)}K_{V_2(z)},h\rangle
	$$
for each $z\in E$ and  $h\in H^2$. Thus
	$$
	W^*K_z=\ov{U_1(z)}K_{V_1(z)}+\ov{U_2(z)}K_{V_2(z)}, \qquad z\in E.
	$$
Since $WC_1K_z=C_1W^*K_z$ for $z\in E$, we get
	$$
	U_1\cdot K_{\bar z}\circ V_1+U_2\cdot K_{\bar z}\circ V_2=U_1(z) K_{\ov{V_1(z)}}+U_2(z) K_{\ov{V_2(z)}}
	$$
	and hence
	\begin{equation}\label{e-new-5}
		\fr{U_1(w)}{1-zV_1(w)}+\fr{U_2(w)}{1-zV_2(w)}=\fr{U_1(z)}{1-V_1(z)w}+\fr{U_2(z)}{1-V_2(z)w}
	\end{equation}
	for all $z\in E$ and $w\in\D$.
By direct computations, we see that (\ref{e-new-5}) yields the following two equalities:
	\begin{equation}\label{e-new-1-1}
		u_1\bar bq(a+b)=(u_2\bar b+v_2\bar a)(\bar a\bar b+qab),
	\end{equation}
	\begin{equation}\label{e-new-1-2}
		u_1q(a+b)+u_1\bar a\bar b^2=(u_2+v_2)qab+(u_2\bar b+v_2\bar a)(\bar a+\bar b).
	\end{equation}
	On the other hand, since $\langle C1, Cz\rangle=0$ and $\langle Cz,1\rangle=\langle C1, z\rangle$, we have
	\begin{equation*}
		u_2+v_2=a(u_2\bar b+v_2\bar a)=u_1\bar a.
	\end{equation*}
Recall $qab=\bar a^2$. It follows from (\ref{e-new-1-1}) and (\ref{e-new-1-2}) that
	\begin{align*}
	\bar b q(a+b)&=a^{-1}\bar a^2(\bar a+\bar b),\\
	q(a+b)&=\bar a(\bar a^2-\bar b^2)+a^{-1}\bar a(\bar a+\bar b).
	\end{align*}
 Therefore, we have $q(a+b)=\bar a^2(\bar a+\bar b)$ and then $a+b=0$ as desired.
	For the case $a=0$ (resp. $b=0$), we can also see  $b=0$ (resp. $a=0$) by the similar argument. We omit the detail
	and finish the proof.
\end{proof}

We close the paper with a characterization of conjugations $C$ on $H^2$ for which
$$
f \in H^\infty \Longrightarrow  CT_fC\in\{T_g: g\in H^\infty\}.
$$
The notation $\mathcal{B}$ denotes the algebra of all bounded linear operators on $H^2$. Also, a Toeplitz operator is called
an analytic Toeplitz operator if its symbol is analytic.

\begin{thm}\label{p-8}
	Let $C$ be a conjugation on $H^2$ and $\Phi$ be a conjugate-linear map on $\mathcal{B}$ defined by
	$\Phi(L)=CLC$ for $L\in \mathcal B$. Then the following statements are equivalent: 
	\begin{enumerate}
\item[(a)]  $\Phi$ leaves the collection of analytic Toeplitz operators invariant.
\item[(b)]  $C=\gamma C_{\lam, a}$ for some $\gamma,\lam\in\T$ and $a\in\D$.
	\end{enumerate}
\end{thm}

\begin{proof}
First assume (b). 
By (\ref{e-new-ee2}), we have
$
		CT_fC=T_{C_\lam(f\circ\va_a)\circ\va_a}$
	for all $f\in H^\infty$ and hence (a) holds.
	
Now assume (a). Then, for each $f\in H^\infty$, there exists $Af\in H^\infty$ such that
	$$
	CT_fC=T_{Af} \qquad {\rm or} \qquad T_f=CT_{Af}C.
	$$ Clearly, such $Af$ is unique.
	This induces a conjugate-linear map $A: H^\infty\to H^\infty$. For each $f\in H^\infty$, note that
	$$
	T_f=CT_{Af}C=C(CT_{AAf}C)C=T_{AAf}.
	$$
Hence $A^2f=f$ for all $f\in H^2$ and then $A^2=I$. Also, for $f,g\in H^\infty$, we note
	\begin{align*}
		T_{A(fg)}
		&=CT_{fg}C=CT_fT_gC\\
		&=(CT_fC)(CT_gC)=T_{Af}T_{Ag}=T_{Af\cdot Ag},
	\end{align*}
which deduce that
	\begin{equation}\label{e-e3}
		A(fg)=Af\cdot Ag.
	\end{equation}
	On the other hand,  since $T_fC=CT_{Af}$, we have
	$$
	f\cdot C1=T_fC1=CT_{Af}1=C(Af\cdot 1)=CAf
	$$
and hence $Af=C(f\cdot C1)$ for all $f\in H^\infty$. Hence we get
	$$
	CT_fC=T_{C(f\cdot C1)}
$$
for all $f\in H^\infty$.
Taking $f=z^n$ and using (\ref{e-e3}),  we get
	$$
	T_{[C(z\cdot C1)]^n}=T_{(Az)^n}=T_{Az^n}=CT_{z^n}C=T_{C(z^n\cdot C1)}
	$$
	for all $n=0,1,2,\cdots$. Also, taking $n=1$ above and using Lemma \ref{lem-new}, we may take $C( z\cdot C1)=\beta\va_{b}$ where
	$\beta\in\T$ and $b\in\D$. Therefore we obtain
	$$
	CT_{z^n}C=T_{(\beta\va_b)^n}, \ \ \ n=0,1,2,\cdots.
	$$
	It follows that
	$
	CT_{z^n}C(e_b)=T_{(\beta\va_b)^n}(e_b)
	$
and hence
	\begin{equation}\label{e-e4}
		C(z^n\cdot Ce_b)=(\beta\va_b)^ne_b
	\end{equation}
for all $n=0,1,2,\cdots$.
	Note that $\{(\beta\va_b)^ne_b\}$ is an \onb~ of $H^2$. Then $\{z^n\cdot Ce_b\}$ is also an \onb~ of $H^2$. Thus
	$$
	\int_\T z^n Ce_b\ov{z^m Ce_b}\, d\sigma=\int_\T z^{n-m}|Ce_b|^2\, d\sigma=0
	$$
	for any nonnegative integers $n,m$ with $n\not=m$. Thus we get
	$$
	\int_\T z^{k}(|Ce_b|^2-1)\, d\sigma=0
	$$
	for any integer $k$ because $||Ce_b||=||e_b||=1$. This implies that $|Ce_b|=1$ a.e. on $\T$ and then
	$Ce_b$ is inner. Since $\{z^n Ce_b\}$ is an \onb~ of $H^2$, we have
	$$
	H^2=\ov{\rm span}\{z^n Ce_b: \  n=0,1,2,\cdots\}=Ce_b\cdot H^2,
	$$
	which implies that $Ce_b=\alpha$, a unimodular constant. Thus (\ref{e-e4}) becomes
	$$
	C z^n=\alpha (\beta\va_b)^ne_b, \ \ \ n=0,1,2,\cdots.
	$$
Now Proposition \ref{p-add} shows $C=\gamma C_{\lam,a}$ for some $\lambda, \gamma\in\T$ and $a\in\D$.
	The proof is complete.
\end{proof}

\section*{Acknowledgements}

The first author was supported by NSFC (11771401) and the second author was supported by Basic Science Research Program through the National Research Foundation of Korea(NRF) funded by the Ministry of Education(NRF-2019R1I1A3A01041943).


\end{document}